\documentclass{agtart_a}
\pdfoutput=1

\title{Amenable groups that act on the line}

\author{Dave Witte Morris}
\givenname{Dave Witte}
\surname{Morris}
\address{Department of Mathematics and Computer Science\\
University of Lethbridge\\\newline
Lethbridge\\Alberta, T1K 3M4\\Canada}
\email{Dave.Morris@uleth.ca}
\urladdr{http://people.uleth.ca/~dave.morris/}

\volumenumber{6}
\issuenumber{}
\publicationyear{2006}
\papernumber{87}
\startpage{2509}
\endpage{2518}

\doi{}
\MR{}
\Zbl{}


\received{9 June 2006}
\revised{}
\accepted{1 September 2006}
\published{15 December 2006}
\publishedonline{15 December 2006}
\proposed{}
\seconded{}
\corresponding{}
\editor{CPR}
\version{}

\arxivreference{math.GR/0606232}

\AtBeginDocument{\let\bar\wbar\let\tilde\wtilde\let\hat\what}
\def\co{\mskip0.5mu\colon\thinspace}

\newcommand{\orders}{\mathcal{O}}
\newcommand{\powerset}{\mathcal{P}}

\newcommand{\real}{{\mathord{\mathbb{R}}}}
\newcommand{\integer}{{\mathord{\mathbb{Z}}}}

\DeclareMathOperator{\SL}{SL}

\numberwithin{equation}{section}

\makeatletter
\def\cnewtheorem#1[#2]#3{\newtheorem{#1}{#3}
\expandafter\let\csname c@#1\endcsname\c@introthm}
\def\dnewtheorem#1[#2]#3{\newtheorem{#1}{#3}[section]
\expandafter\let\csname c@#1\endcsname\c@thm}
\makeatother

\newtheorem{introthm}{Theorem}
\makeautorefname{introthm}{Theorem}

\cnewtheorem{introcor}[introthm]{Corollary}
\makeautorefname{introcor}{Corollary}

\theoremstyle{remark}
\newtheorem*{introrems}{Remarks}

\theoremstyle{plain}
\newtheorem{thm}{Theorem}[section]
\dnewtheorem{lem}[thm]{Lemma}
\dnewtheorem{cor}[thm]{Corollary}
\dnewtheorem{prop}[thm]{Proposition}
\dnewtheorem{questions}[thm]{Questions}

\theoremstyle{definition}
\dnewtheorem{defn}[thm]{Definition}
\dnewtheorem{eg}[thm]{Example}

\theoremstyle{remark}
\newtheorem*{ack}{Acknowledgments\ignorespaces}
\dnewtheorem{rem}[thm]{Remark}

\begin{document}

\begin{asciiabstract}
Let Gamma be a finitely generated, amenable group.  Using an idea of E
Ghys, we prove that if Gamma has a nontrivial, orientation-preserving
action on the real line, then Gamma has an infinite, cyclic quotient.
(The converse is obvious.)  This implies that if Gamma has a faithful
action on the circle, then some finite-index subgroup of Gamma has the
property that all of its nontrivial, finitely generated subgroups have
infinite, cyclic quotients.  It also means that every left-orderable,
amenable group is locally indicable.  This answers a question of P
Linnell.
\end{asciiabstract}

\begin{htmlabstract}
Let &Gamma; be a finitely generated, amenable group.  Using an idea
of &Eacute;&nbsp;Ghys, we prove that if &Gamma; has a nontrivial,
orientation-preserving action on the real line, then &Gamma; has an
infinite, cyclic quotient.  (The converse is obvious.)  This implies
that if &Gamma; has a faithful action on the circle, then some
finite-index subgroup of&nbsp;&Gamma; has the property that all of its
nontrivial, finitely generated subgroups have infinite, cyclic
quotients.  It also means that every left-orderable, amenable group is
locally indicable.  This answers a question of P&nbsp;Linnell.
\end{htmlabstract}

\begin{abstract} 
Let $\Gamma$ be a finitely generated, amenable group.  Using an idea
of \'E~Ghys, we prove that if $\Gamma$ has a nontrivial,
orientation-preserving action on the real line, then $\Gamma$ has an
infinite, cyclic quotient.  (The converse is obvious.)  This implies
that if $\Gamma$ has a faithful action on the circle, then some
finite-index subgroup of~$\Gamma$ has the property that all of its
nontrivial, finitely generated subgroups have infinite, cyclic
quotients.  It also means that every left-orderable, amenable group is
locally indicable.  This answers a question of P~Linnell.
\end{abstract}

\maketitle

\setcounter{section}{-1}

\section{Introduction}

Let $\Gamma$ be an abstract group (with the discrete topology). It is
obvious that if $\Gamma$ has an infinite cyclic quotient, then
$\Gamma$ has a nontrivial, orientation-preserving action on the real
line~$\real$.  The converse is not true in general, even for finitely
generated groups \cite[Example 6.9.2, page 128]{GlassPOG}.  In
this note, we use an idea of \'E~Ghys to prove that the converse does
hold in the class of finitely generated, amenable groups.

\begin{defn}\cite[page 9 and Theorem 5.4(i,iii), page 45]{Pier}\label{AmenableDefn} \ 
 \begin{itemize}
\item A measure~$\mu$ on a measure space~$X$ is said to be a
\emph{probability measure} iff $\mu(X) = 1$.
\item A (discrete) group~$\Gamma$ is \emph{amenable} iff for every continuous
action of~$\Gamma$ on a compact, Hausdorff space~$X$, there is a
$\Gamma$--invariant probability measure on~$X$.
\end{itemize}
 \end{defn}

\begin{introthm} \label{ActR<>CyclicQuot}
Let\/ $\Gamma$ be a finitely generated, amenable group. Then $\Gamma$ has a nontrivial, orientation-preserving action on\/~$\real$ if and only if\/ $\Gamma$ has an infinite cyclic quotient.
\end{introthm}

It is well known that a countable group has a faithful,
 orientation-preserving action on~$\real$ if and only if it is
 \emph{left orderable} \cite[Theorem 6.8]{Ghys-AOC}.  (That is,
 there is a \emph{left-invariant order} on~$\Gamma$; in other words,
 there is a total order~$\prec$ on~$\Gamma$, such that, for all
 $\gamma, \lambda_1,\lambda_2 \in \Gamma$, if $\lambda_1 \prec
 \lambda_2$, then $\gamma \lambda_1 \prec \gamma \lambda_2$.)  Also,
 every subgroup of an amenable group is amenable \cite[Proposition
 13.3]{Pier}.  Hence, the nontrivial direction of
 \fullref{ActR<>CyclicQuot} can be stated in the following purely
 algebraic form.

\begin{defn}{\cite[page 127]{GlassPOG}}\qua
A group is \emph{locally indicable} iff each of its nontrivial finitely generated subgroups has an infinite cyclic quotient.
\end{defn}

\begin{introthm} \label{Amenable->LocInd}
Every amenable left-orderable group is locally indicable.
\end{introthm}

\begin{introrems} \ 
\begin{enumerate}
\item  The theorem answers  a question of P~Linnell \cite[page 134]{Linnell-LOAmen}.
\item Every locally indicable group (whether amenable or not) is left
  orderable (Burns and Hale \cite{BurnsHale}, \cite[Lemma
  6.9.1]{GlassPOG}).
\item \fullref{Amenable->LocInd} has previously been proved with
``amenable" replaced by stronger hypotheses, such as polycyclic
(Rhemtulla \cite{Rhemtulla-polycyclic}), solvable-by-finite (Chiswell
and Kropholler \cite{ChiswellKropholler}), supramenable (Kropholler
\cite{Kropholler-AmenRO}), or elementary amenable (Linnell
\cite{Linnell-LOAmen}).  There are also interesting results that
replace ``amenable" with the assumption that $\Gamma$ has no
nonabelian free \emph{semi}groups
(Longobardi, Maj and Rhemtulla \cite{LongobardiMajRhemtulla-nofree}), or generalizations of this
(Longobardi, Maj and Rhemtulla \cite{LongobardiMajRhemtulla-when},
Linnell \cite{Linnell-LOnofree}).
\item P~Linnell \cite[Conjecture 1.1]{Linnell-LOnofree} has conjectured that the theorem is valid with ``amenable" replaced by the weaker condition of not containing any nonabelian free subgroups.
\end{enumerate}
\end{introrems}

Because the universal cover of the circle~$S^1$ is~$\real$, the following corollary is an easy consequence (cf \cite[Section 5]{Linnell-LOnofree}).

\begin{introcor} \label{ActCircle->VirtLocind}
If\/ $\Gamma$ is an amenable group, and\/ $\Gamma$ has a faithful, orientation-preserving action on~$S^1$, then $\Gamma$ has a normal subgroup~$N$, such that $N$ is locally indicable, and $\Gamma/N$ is a finite cyclic group.
\end{introcor}

\begin{proof}[Proof of \fullref{Amenable->LocInd}]
 Let $\Gamma$ be an amenable group that is left orderable. 
We wish to show that $\Gamma$ is locally indicable,
so there is no harm in assuming that $\Gamma$ is nontrivial
and finitely generated (hence, countable).
Details of each of the following steps of the proof (and
the necessary definitions) are presented in the indicated section
below.
 \begin{enumerate}\leftskip 21pt
 \item[\fullref{SpaceOfOrdersSect}] Let $\orders$ be the collection of all
left-invariant orders on~$\Gamma$. \'Etienne Ghys observed that $\orders$ 
is a compact Hausdorff
space, and that the action of~$\Gamma$ on~$\orders$ by right translations is
continuous.
 \item[\fullref{AmenableSect}] Since $\Gamma$ is amenable, there is a
$\Gamma$--invariant probability measure~$\mu$ on~$\orders$.
 \item[\fullref{PoincareSect}] The Poincar\'e Recurrence Theorem implies
there is a point in~$\orders$ that is recurrent for every cyclic subgroup
of~$\Gamma$.
 \item[\fullref{LocInd&RecSect}] Being recurrent for every cyclic subgroup is a
stronger condition than being Conradian, and it is well known that any
group admitting a Conradian left-invariant order is locally indicable.
 \end{enumerate}
 Therefore, $\Gamma$ is locally indicable.
 \end{proof}

\begin{ack}
 I would like to thank \'Etienne Ghys for pointing out (several years ago)
that the space of orderings of~$\Gamma$ is compact, and for suggesting
that it would be worthwhile to study the action of~$\Gamma$ on this space.
 I also benefitted from conversations with Uri Bader that clarified and
extended Ghys' observations, and from further discussions with Alex
Furman and Tsachik Gelander. 
Peter Linnell provided helpful comments on a preliminary version of this manuscript.

 I would also like to thank the Department of Mathematics at the
University of Chicago for their hospitality while this research was being
carried out.
 The work was partially supported by a grant from the National Science and
Engineering Research Council of Canada.
 \end{ack}

\section{The space of left-invariant orders}
 \label{SpaceOfOrdersSect}

 Fix a left-orderable group~$\Gamma$.
 In this section, we present an idea of \'E~Ghys (personal communication).
 The space of left-invariant orderings was topologized in a paper of A\,S~Sikora \cite{Sikora}, but the action of~$\Gamma$ on this space does not seem to have appeared previously in the literature.

\begin{defn} \ 
 \begin{enumerate}
 \item Let $\orders$ be the collection
of all left-invariant orders on~$\Gamma$. 
 (Note that $\orders$ is nonempty, because $\Gamma$ is left orderable.)
 \item (Ghys, Sikora \cite[Definition 1.1]{Sikora})\qua
Topologize $\orders$ by declaring a set to be open iff it is a
union of basic open sets of the form
 $$ U_{\lambda_1,\lambda_2,\ldots,\lambda_r} = \{\, {\prec} \mid
\lambda_1 \prec \lambda_2 \prec \cdots \prec \lambda_r \,\} $$
 for a sequence of distinct $\lambda_1,\lambda_2,\ldots,\lambda_r \in
\Gamma$.
 \item (Ghys)\qua Let $\Gamma$ act on~$\orders$ (on the right) by right
translation:
	$$ \lambda_1 \prec_\gamma \lambda_2 \iff \lambda_1 \gamma^{-1} \prec
\lambda_2 \gamma^{-1} .$$
 It is clear that if $\prec$ is a left-invariant order, then
$\prec_\gamma$ is a left-invariant order, for every $\gamma \in \Gamma$.
Also, we have ${\prec_{\gamma_1 \gamma_2}} =
({\prec_{\gamma_1}})_{\gamma_2}$, so this defines an action of~$\Gamma$.
 \end{enumerate}
 \end{defn}

\begin{lem} \label{OrderSpaceLemma} \ 
 \begin{enumerate}
 \item \label{OrderSpaceLemma-cpct}
 \text{\upshape(Ghys, Sikora)}\qua
  $\orders$ is a compact, Hausdorff space.
 \item \label{OrderSpaceLemma-homeo}
 \text{\upshape(Ghys)}\qua
 $\Gamma$ acts on~$\orders$ by homeomorphisms.
 \end{enumerate}
 \end{lem}

\begin{proof}

 \ref{OrderSpaceLemma-cpct}\qua \cite[Theorem 1.4]{Sikora}\qua
 The collection $\powerset(\Gamma \times \Gamma)$
of all subsets of $\Gamma \times \Gamma$ is a compact, Hausdorff space
(because it is naturally homeomorphic to the infinite Cartesian product
$\{0,1\}^{\Gamma \times \Gamma}$).
It is easy to see that the complement of~$\orders$ is open, so $\orders$ is compact.

\ref{OrderSpaceLemma-homeo}\qua The image of the basic open set
$U_{\lambda_1,\lambda_2,\ldots,\lambda_r}$ under an element~$\gamma$
of~$\Gamma$ is the basic open set
 $U_{\lambda_1 \gamma,\lambda_2 \gamma,\ldots,\lambda_r \gamma}$.
 \end{proof}

\section{Amenability}
 \label{AmenableSect}

The following observation is immediate from \fullref{OrderSpaceLemma} and the definition of amenability.

\begin{lem} \label{Amenable->MeasOnO}
 Let $\Gamma$ be a left-orderable, amenable group. Then there
is a $\Gamma$--invariant probability measure on the space~$\orders$ of
left-invariant orders on~$\Gamma$.
 \end{lem}

\begin{rem}
 The above lemma is the only use that will be made of amenability. Hence,
in the statements of \fullref{ActR<>CyclicQuot}, \fullref{Amenable->LocInd}, 
and \fullref{ActCircle->VirtLocind}, amenability can be
replaced with the assumption that there is a $\Gamma$--invariant
probability measure on~$\orders$.
 \end{rem}

\section{Poincar\'e Recurrence Theorem}
 \label{PoincareSect}

We recall the following classical result that can be
found in almost any textbook on Ergodic Theory.
For the convenience of the reader, we include a short proof.

\begin{prop}{\rm(Poincar\'e Recurrence Theorem \cite[page 7]{SinaiTopics})}\qua
 Suppose \begin{itemize} \item $X$ is a measure space with probability
 measure~$\mu$, \item $T \co X \to X$ is an invertible, measurable map
 that preserves the measure~$\mu$, and \item $A$ is any measurable
 subset of~$X$.  \end{itemize} Then there is a measurable subset~$Z$
 of~$X$ with $\mu(Z) = 0$, such that, for every $a \in A
 \smallsetminus Z$, there is a sequence of positive integers $n_i \to
 \infty$, such that $T^{n_i}(a) \in A$ for every~$i$.  \end{prop}
 
\begin{proof}
For $n \in \integer^+$, let 
	$$ A_n = \bigcup_{k=1}^\infty T^{-k n}(A) .$$
For each $a \in \bigcap_{n=1}^\infty A_n$, there is a sequence of
positive integers $n_i \to \infty$, such that $T^{n_i}(a) \in A$ for
every~$i$. Thus, it suffices to show that $\mu(A \smallsetminus A_n) = 0$.

Suppose $\mu(A \smallsetminus A_n) > 0$. Since the sets 
	$$T^{-n}(A \smallsetminus A_n), \ T^{-2n}(A \smallsetminus A_n), \ T^{-3n}(A \smallsetminus A_n), \ \ldots$$
all have the same measure (and $\mu(X) < \infty$), they cannot all be disjoint. Hence, there exist $k > \ell$, such that $T^{-kn}(A \smallsetminus A_n) \cap T^{-\ell n}(A \smallsetminus A_n) \neq \emptyset$. By applying $T^{\ell n}$, we may assume $\ell = 0$. Therefore
 	$$\emptyset 
\neq T^{-kn}(A \smallsetminus A_n) \cap (A \smallsetminus A_n) 
\subset T^{-kn}(A) \cap \bigl( A \smallsetminus T^{-kn}(A) \bigr) 
=\emptyset
.$$
This is a contradiction.
\end{proof}

\begin{defn} \label{RecurrentOrderDefn}
 A left-invariant order~$\prec$ on a group~$\Gamma$ is
\emph{recurrent for every cyclic subgroup} iff for every $\gamma \in \Gamma$ and
every (finite) increasing sequence 
$\lambda_1 \prec \lambda_2 \prec \cdots \prec \lambda_r$ of elements of~$\Gamma$, 
there exist positive integers $n_i \to
\infty$, such that
 $$ \lambda_1 \gamma^{n_i} \prec \lambda_2 \gamma^{n_i} \prec \cdots
\prec \lambda_r \gamma^{n_i} $$
 for every~$i$.
 \end{defn}

\begin{cor}
 Suppose
 \begin{itemize}
 \item $\Gamma$ is a left-orderable group that is countable, 
 and
 \item there exists a $\Gamma$--invariant probability measure~$\mu$ on
the space~$\orders$ of left-invariant orders on~$\Gamma$.
 \end{itemize}
 Then $\Gamma$ admits a left-invariant order that is recurrent for every cyclic subgroup.
 \end{cor}

\begin{proof}
 For each $\gamma \in \Gamma$ and each sequence
$\lambda_1,\lambda_2,\ldots,\lambda_r$ of distinct elements of~$\Gamma$,
we may apply the Poincar\'e Recurrence Theorem with
 \begin{itemize}
 \item the space~$\orders$ in the role of~$X$,
 \item the transformation~$\gamma^{-1}$ in the role of~$T$,
 and
 \item the basic open set $U_{\lambda_1,\ldots,\lambda_r}$ in the role
of~$A$. 
 \end{itemize}
 We conclude that there is a
set~$Z_{\gamma,\lambda_1,\ldots,\lambda_r}$ of measure~$0$
in~$\orders$, such that:
\begin{quote}
 if $\prec$ is any left-invariant order
on~$\Gamma$, such that
 \begin{itemize}
 \item $\lambda_1 \prec \lambda_2 \prec \cdots \prec \lambda_r$,
 and
 \item $\prec$ is not in~$Z_{\gamma,\lambda_1,\ldots,\lambda_r}$,
 \end{itemize}
 then there exist positive integers $n_i \to \infty$, such that
 $$ \text{$\lambda_1 \gamma^{n_i} \prec \lambda_2 \gamma^{n_i} \prec
\cdots \prec \lambda_r \gamma^{n_i}$ \quad for every~$i$.}	$$ 
\end{quote}
 The union of all of the sets $Z_{\gamma,\lambda_1,\ldots,\lambda_r}$ has
measure zero (because it is a countable union of sets of measure~$0$), so
there is a left-invariant order~$\prec$ that does not belong to
any~$Z_{\gamma,\lambda_1,\ldots,\lambda_r}$. This order is recurrent for every cyclic subgroup.
 \end{proof}

Combining the above corollary with \fullref{Amenable->MeasOnO}
immediately yields the following conclusion.

\begin{cor}
 If\/ $\Gamma$ is any countable, left-orderable, amenable group, 
 then $\Gamma$ admits a left-invariant order that is recurrent for every cyclic subgroup.
 \end{cor}

\section{Recurrent orders and indicable groups}
 \label{LocInd&RecSect}
 
The main result of this section is \fullref{Recurrent->LocInd}. It is an almost immediate 
consequence of a known result (\fullref{LocInd<>Conradian} below),
but, for the convenience of the reader, we provide a short proof that is fairly self contained.

\begin{defn}\cite[Lemma 6.6.2(1,3), page 121]{GlassPOG}\qua
 A left-invariant order~$\prec$ on a group~$\Gamma$ is
\emph{Conradian} iff for every $\gamma,\lambda \in \Gamma$, such that
$\gamma \succ e$ and $\lambda \succ e$, there exists $n \in \integer^+$,
such that $\lambda \gamma ^n \succ \gamma $.
 \end{defn}

\begin{thm}{\rm\cite[Theorem 6.J, page 128]{GlassPOG}}\qua
\label{LocInd<>Conradian}
 A group is locally indicable if and only if it admits a Conradian
left-invariant order.
 \end{thm}

\begin{lem} \label{Recurrent->Conradian}
 If a left-invariant order is recurrent for every cyclic subgroup, 
 then the order is Conradian.
 \end{lem}

\begin{proof}
 If $\lambda \succ e$, then recurrence implies there exists $n \in
\integer^+$, such that $\lambda \gamma ^n \succ e \gamma ^n = \gamma ^n$.
If, in addition, we have $\gamma \succ e$, then $\gamma ^n \succeq
\gamma $, so transitivity implies $\lambda \gamma ^n \succ \gamma $.
 \end{proof}

Combining \fullref{Recurrent->Conradian} with
\fullref{LocInd<>Conradian} yields the following result:
 
\begin{cor} \label{Recurrent->LocInd}
 If a group admits a left-invariant order that is recurrent for every cyclic subgroup,
 then the group is locally indicable.
 \end{cor}

Before providing a proof of 
\fullref{Recurrent->LocInd} that does not rely on \fullref{LocInd<>Conradian},
let us recall some elementary properties of left-ordered groups.

\begin{rem} 
Let $\prec$ be a left-invariant order on a group~$\Gamma$.
\begin{enumerate}
\item To say that a subgroup~$\Lambda$ of~$\Gamma$ is \emph{convex}
means that if $\lambda_1 \prec \gamma \prec \lambda_2$, with 
$\lambda_1,\lambda_2 \in \Lambda$, then $\gamma \in \Lambda$ \cite[page 31]{GlassPOG}.
(Because $\prec$ is left invariant, it suffices to verify this condition in the special case
where $\lambda_1 = e$.)
\item  If $\Lambda$ is a convex, proper subgroup of~$\Gamma$, then:
\begin{enumerate}
\item $\Lambda$ is an interval in the total order $(\Gamma,\prec)$ (so left invariance implies that each left coset of~$\Lambda$ is also an interval),
so
\item $\prec$ induces a well-defined total order on the space of left cosets of~$\Lambda$,
so
\item $\Lambda$ is bounded above (by any positive element of~$\Gamma$ that does not belong to~$\Lambda$) and bounded below.
\end{enumerate}
\item If $\Gamma$ is finitely generated, then Zorn's Lemma implies that $\Gamma$ has a maximal (proper) convex subgroup. The convex subgroups of~$\Gamma$ are totally ordered by inclusion \cite[Lemma 3.1.2, page 32]{GlassPOG}, so this maximal convex subgroup is unique.
\item The order $\prec$ is \emph{Archimedean} if, for all nontrivial $\gamma,\lambda \in \Gamma$, 
there exists $n \in \integer$, such that $\gamma \preceq \lambda^n$ \cite[page 55]{GlassPOG}. 
It is well known (and not difficult to show \cite[Corollary 4.1.3, page 56]{GlassPOG}) that any group admitting an Archimedean left-invariant order must be abelian (and torsion free). 
\end{enumerate}
\end{rem}
 
 \begin{proof}[Direct proof of \fullref{Recurrent->LocInd}]
Assume that $\prec$ is a left-invariant order on a finitely generated group~$\Gamma$,
and that $\prec$ is recurrent for every cyclic subgroup.

Begin by noting, for $\lambda_1,\lambda_2, \ldots, \lambda_r, \lambda \in \Gamma$, 
that if $\max \bigl\{ \lambda_1^{\pm1}, \ldots, \lambda_r^{\pm1} \bigr\} \preceq \lambda$, then there exists
$n \in \integer^+$, such that 
	$$\lambda_1 \lambda_2 \cdots \lambda_r \preceq \lambda^n .$$
To see this, choose (by induction on~$r$) some $m \in \integer^+$, such that 
$\lambda_2 \lambda_3 \cdots \lambda_r \preceq \lambda^m$, and then choose 
(by recurrence of~$\prec$) some $n \in \integer^+$, such that 
	$$\lambda_1 \lambda^{mn} \preceq \lambda^m \lambda^{mn} = \lambda^{m(n+1)} .$$
Then we have $\lambda_2 \lambda_3 \cdots \lambda_r \preceq \lambda^m \preceq \lambda^{mn}$, so
	$$ \lambda_1 \lambda_2 \lambda_3 \cdots \lambda_r \preceq \lambda_1 \lambda^{mn} \preceq \lambda^{m(n+1)} ,$$
as desired.

Since $\Gamma$ is finitely generated, it has a maximal (proper) convex subgroup~$\Lambda$,
which is unique. For any $\gamma \in \Gamma^+$, the set
	$$ \Lambda_\gamma = \{\, \lambda \in \Gamma \mid \text{$\lambda^n \prec \gamma$ for all $n \in \integer$} \,\} $$
is obviously closed under inverses, and the observation of the preceding paragraph implies that it
is closed under multiplication and that it is convex; hence, it follows from the uniqueness of the maximal convex subgroup that $\Lambda_\gamma \subset \Lambda$.
Thus, 
 \begin{equation} \tag{$*$} \label{BddAboveInLambda}
 \text{$\Lambda$ contains every subgroup of~$\Gamma$ that is bounded above.}
 \end{equation}
The reversal of a recurrent order is recurrent, so the same argument implies
 \begin{equation} \tag{$**$} \label{BddBelowInLambda}
 \text{$\Lambda$ contains every subgroup of~$\Gamma$ that is bounded below.}
 \end{equation}
Since $\Lambda$ is a convex proper subgroup, it is bounded both above (by some~$\lambda_+$) and below (by some~$\lambda_-$). For $\gamma \in \Gamma$, and any $\lambda \in \Lambda$, we have
	$$ \gamma \succ e 
	\implies \gamma^{-1} \prec e
	\implies \lambda \gamma^{-1} \prec \lambda \preceq \lambda_+
	\implies \gamma \lambda \gamma^{-1} \prec \gamma \lambda_+ .$$
Similarly,
	$$ \gamma \prec e \implies \gamma \lambda \gamma^{-1} \succ \gamma \lambda_- .$$
Thus, for any $\gamma \in \Gamma$, the conjugate $\gamma \Lambda \gamma^{-1}$ is either bounded above (by~$\gamma \lambda_+$) or bounded below (by~$\gamma \lambda_-$). From \ref{BddAboveInLambda} and \ref{BddBelowInLambda}, we conclude that $\gamma \Lambda \gamma^{-1} \subset \Lambda$. Therefore, $\Lambda$ is normal in~$\Gamma$.

The order induced on the quotient group $\Gamma/\Lambda$ is Archimedean (because \ref{BddAboveInLambda} implies no nontrivial subgroup of $\Gamma/\Lambda$ is bounded above), so $\Gamma/\Lambda$ must be abelian and torsion free. Then the structure of finitely generated abelian groups implies that $\Gamma/\Lambda$ has an infinite cyclic quotient.
\end{proof}

 The converse of \fullref{Recurrent->LocInd} is false:

\begin{eg}
 Let $F$ be a free subgroup of finite index in $\SL(2,\integer)$,
and let $\Gamma$ be the natural semidirect product $F \ltimes
\integer^2$. Then $\Gamma$ is locally indicable, but has no left-invariant
order that is recurrent for every cyclic subgroup.
 \end{eg}
 
 \begin{proof}
 Since free groups and $\integer^2$ are locally indicable, it is clear that
$\Gamma$ is locally indicable.

Suppose $\prec$ is a left-invariant order on~$\Gamma$
that is recurrent for every cyclic subgroup.
(This will lead to a contradiction.)
For a matrix~$T$ in~$F$ and a vector~$v$ in~$\integer^2$, 
let us use $\overline{T}$ and $\overline{v}$ to represent the 
corresponding elements of $F \ltimes \integer^2$, so
$\overline{T} \, \overline{v} \, \overline{T}^{-1} 
= \overline{T(v)}$.

There is a nonzero linear functional $f \co \real^2 \to \real$,
such that $f(v) > 0 \Rightarrow \overline{v} \succ e$, for all 
$v \in \integer^2$.
Let $T$ be a hyperbolic matrix in~$F$, with eigenvalues
$\alpha_1,\alpha_2>0$ and corresponding eigenvectors $v_1$ and~$v_2$.
Since $v_1$ and~$v_2$ are linearly independent,
we may assume $f(v_1) \neq 0$. Furthermore, we may assume 
$f(v_1) > 0$ and $\alpha_1 > 1$ (so $\alpha_2 < 1$), by replacing 
$v_1$ with~$-v_1$
and/or $T$ with~$T^{-1}$, if necessary. Let $g \co \real^2 \to \real$
be the (unique) linear functional that satisfies $g(v_1) = 1$ and $g(v_2) = 0$.

Given any $v \in \integer^2$, such that $f(v) > 0$, recurrence
provides a sequence $n_i \to \infty$, such that 
$\overline{v} \, \overline{T}^{-n_i} \succ \overline{T}^{-n_i}$ 
for every~$i$. By left invariance, this implies 
$\overline{T^{n_i}(v)} 
= \overline{T}^{n_i} \overline{v} \overline{T}^{-n_i} \succ e$,
so 
$ f \bigl( T^{n_i}(v) \bigr) \ge 0 $.
Since
$$ \lim_{n \to \infty} \frac{1}{\alpha_1^n} f \bigl( T^{n}(v) \bigr) 
 = f \left( \lim_{n \to \infty} \frac{1}{\alpha_1^n} T^{n}(v) \right)
  = f \bigl( g(v) v_1 \bigr)
   = g(v) \, f(v_1) ,$$
we conclude that $g(v) \ge 0$. 

Since $v$ is an arbitrary element of~$\integer^2$ with $f(v) > 0$, 
this implies that $\ker f = \ker g$ is an eigenspace of~$T$. 
But $T$ is an arbitrary hyperbolic matrix in~$F$, and it is easy to
show that there are two (conjugate) hyperbolic matrices in~$F$
that do not have a common eigenspace. This is a contradiction.
\end{proof}

\bibliographystyle{gtart}
\bibliography{link}

\end{document}